\documentclass[a4paper]{article}

\usepackage{amsmath}
\usepackage{amssymb}

\newtheorem{thm}{Theorem}

\newtheorem{prp}[thm]{Proposition}

\newtheorem{lem}[thm]{Lemma}

\newcommand{\nn}{\nonumber}
\newcommand{\p}{\partial}
\newcommand{\vp}{\varphi}
\newcommand{\ve}{\epsilon}
\newcommand{\dl}{\delta}

\newcommand{\pal}{\partial}

\newcommand{\al}{\alpha}

\newenvironment{prf}{\noindent {\it Proof} \ }{\hfill $\Box$}
\newenvironment{prfM}{\noindent {\it Proof of the Main Theorem} \ }{\hfill $\Box$}
\begin{document}
\title {On the Quasitriviality of Deformations of Bihamiltonian Structures
of Hydrodynamic Type}

\author{{Si-Qi Liu \ \ Youjin Zhang}\\
{\small Department of Mathematical Sciences, Tsinghua
University}\\
{\small Beijing 100084, P.R.China}\\
{\small  lsq99@mails.tsinghua.edu.cn,\
yzhang@math.tsinghua.edu.cn}}
\date{}
\maketitle

\begin{abstract}
We prove in this paper the quasitriviality of a class of deformations of the one component
bihamiltonian structures of hydrodynamic type.
\end{abstract}

\section{Introduction}
The notion of quasitrivial deformations of a bihamiltonian structure of hydrodynamic type is
introduced in \cite{DZ1}, it is proved there that for certain class of deformations of a semisimple
bihamiltonian structure of hydrodynamic type\cite{dn83, dn84, dn89, maltsev} there exist quasi-Miura transformations which
transform the deformed bihamiltonian structures to the undeformed ones. The proof of quasitriviality
that is given in \cite{DZ1} uses the requirement that
the deformed bihamiltonian structures satisfy the so call tau symmetry condition, roughly speaking, it
requires that the hierarchy of integrable systems that are related to the deformed bihamiltonian
structures possess tau functions. The purpose of this paper is to show, without the assumption of
tau symmetry, that the deformations of the
following bihamiltonian structure of hydrodynamic type
\begin{eqnarray}
&&\{u(x),u(y)\}_1=\varphi(u(x))\,\delta'(x-y)+\frac12\, \p_x(\vp(u(x))\,\delta(x-y),\nn\\
&&\{u(x),u(y)\}_2=u(x)\,\vp(u(x))\,\delta'(x-y)+\frac12\,{\p_x(u(x) \vp(u(x))}\,\delta(x-y) \label{z2}
\end{eqnarray}
are quasitrivial. Here $\vp(u)$ is an arbitrary smooth function, and it is easy to see that this is the
general form of a one component bihamiltonian structure of hydrodynamic type.
The deformations that we are to study are assumed to take the form
\begin{eqnarray}
&&\{u(x),u(y)\}_{a}=\{u(x),u(y)\}^{[0]}_{a}
\nn\\&&+\sum_{m\ge 1}\sum_{l=0}^{m+1} \ve^m A_{m,l;a}(u;u_x,\dots,u^{(m+1-l)}) \delta^{(l)}(x-y),
\quad  a=1,2.\label{z1}
\end{eqnarray}
Here $\{\,,\,\}^{[0]}_a$ denote the undeformed Poisson brackets (\ref{z2}), the functions $A_{m,l;a}$
are polynomials in $u_x, u_{xx},\dots, u^{(m+1-l)}=\p^{m+1-l}_x u(x)$, and the coefficients
of these polynomials are smooth functions of $u(x)$. The functions $A_{m,l;a}$ are also required to be
homogeneous of degree $m+1-l$, here we assign a degree
$k$ to $u^{(k)}=\p^k_x u(x)$.
When the function $\vp(u)$ is taken to be a constant, the above deformed bihamiltonian structures include
the well known examples that correspond to the KdV hierarchy
and the Camassa-Holm hierarchy \cite{CH,CHH, Fu, gardner, magri, ZF}, and they also contain some more
general examples of deformations that are considered in \cite{Lo}.
For example, let us take $\vp(u)=1$, then we have the following bihamiltonian structure
of the KdV hierachy
\begin{eqnarray}
&&\{u(x),u(y)\}_1=\delta'(x-y),\nn\\
&&\{u(x),u(y)\}_2=u(x)\delta'(x-y)+\frac12\,u(x)'\delta(x-y)+\frac{\ve^2}{8}\,\delta'''(x-y),\label{z2b}
\end{eqnarray}
The quasitriviality of this bihamiltonian structure means that under the quasi-Miura transformation
\begin{equation}
u\mapsto u +{\epsilon^2\over 24} \p_x^2 \left( \log u_x\right)
+\epsilon^4 \p_x^2 \left( {u^{(4)}\over 1152\, {u'}^2}
- {7\, u'' u'''\over 1920\, {u'}^3}
+{{u''}^3\over 360\, {u'}^4}\right) + O(\epsilon^6)
\end{equation}
it is transformed to the bihamiltonian structure (\ref{z2b}) with $\ve=0$ \cite{DZ1}.

The main result of the paper is the following\newline
\vskip -0.2truecm
\noindent {\bf Main Theorem.}\ {\em
Any deformation (\ref{z1}) of the bihamiltonian structure (\ref{z2}) can be obtained from the undeformed one by
a quasi-Miura transformation of the form
\begin{equation}\label{qquasi-Miura}
u\mapsto u+\sum_{k\ge 1} \ve^k G_k(u;u_x,\dots,u^{(m_k)}).
\end{equation}
Here $G_k$ are smooth functions of their arguments.}

If a deformation of the bihamiltonian structure (\ref{z2}) is obtained from a quasi-Miura transformation (\ref{qquasi-Miura})
with coefficients $G_k, k\ge 1$ that are polynomials in $u_x,\dots, u^{(m_k)}$, then such a deformation is called
trivial \cite{DZ1}, and the related quasi-Miura transformation is called a Miura transformation.

We will prove the Main Theorem in section 2, and in section 3 we will discuss generalizations of the above result
to the case of multi-component bihamiltonian structures.

\section{Proof of the Main Theorem}
Let us  first recall the definition of the space of local
multi-vectors and the operation of Schouten-Nijenhuis bracket that is defined on it
(see for details in \cite{DZ1} and references therein ).
A local $k$ vector is defined to be a formal infinite sum of the following form
\begin{equation}
\al=\sum \frac1{k!}\pal_{x_1}^{s_1}\dots\pal_{x_k}^{s_k} A\,
\frac{\pal}{\pal u^{(s_1)}(x_1)}\wedge\dots \wedge\frac{\pal}{\pal u^{(s_k)}(x_k)}
\end{equation}
with the coefficient $A$ having the expression
\begin{equation}
A=\sum_{p_2,\dots,p_{k}\ge 0} B_{p_2\dots p_k}(u(x_1);
u_x(x_1),\dots) \delta^{(p_2)}(x_1-x_2)\dots \delta^{(p_k)}(x_1-x_k).
\end{equation}
Here $B_{p_2\dots p_k}$ are smooth functions of $u, u_x,\dots, u^{(m)}$ for certain integers
$m$ which may depend on the indices $p_2,\dots,p_k$, and
\begin{equation}
A=A(x_1,\dots,x_k;u(x_1),\dots,u(x_k),\dots)
\end{equation}
is antisymmetric with respect to the permutations
$
x_p\leftrightarrow x_q.
$
The distribution  $A$ is called the component of the local $k$-vector $\al$.
We denote by $\Lambda_{loc}^k$ the space of all such local $k$-vectors. {\em Note that in the
definition of a local $k$-vector that is given in
\cite{DZ1} the functions $B_{p_2\dots p_k}$ are required to be differential polynomials}.
We drop here the polynomiality condition for the convenience of later use, nevertheless, we
still keep to use the notations of \cite{DZ1} such as $\Lambda^k_{loc}$.
In particular, a local vector field has the expression
\begin{equation}\label{def-vec}
\xi=\sum_{s\ge 0} \pal_x^s X(u(x);u_x(x),\dots,u^{(m)})\frac{\pal}{\pal u^{(s)}(x)}
\end{equation}
and a local bivector takes the form
\begin{equation}\label{bi-vec}
\omega=\frac12\sum \pal_{x}^{s} \pal_y^t A\frac{\pal}{\pal u^{(s)}(x)}\wedge\frac{\pal}{\pal u^{(t)}(y)}
\end{equation}
with
\begin{equation}
A=\sum_{k\ge 0} A_k(u(x);u_x(x),\dots,u^{(m_k)}) \delta^{(k)}(x-y).
\end{equation}
We also denote by $\Lambda_{loc}^0$ the space that consists of
local functionals of the form
\begin{equation}\label{def-func}
{\bar f}=\int f(u(x);u_x(x),\dots,u^{(m)}) dx.
\end{equation}

The operation of Schouten-Nijenhuis bracket
\begin{equation}
[\ ,\, ]: \ \Lambda^k_{loc}\times \Lambda^{l}_{loc}\to\Lambda^{k+l-1}_{loc},\quad k,l\ge 0
\end{equation}
generalizes that of the commutators between vector fields. For example,
the components of the Schouten-Nijenhuis bracket of a bivector $\omega$ of the form (\ref{bi-vec})
with a local functional $\bar f$ and with a local vector filed $\xi$ of the form (\ref{def-vec})
are given respectively by
\begin{equation}
\sum_{k} A_k \pal_x^k\frac{\delta \bar f}{\delta u(x)}
\end{equation}
and
\begin{eqnarray}
&&\sum_{t}\left(\pal_x^t X(u(x); \dots)
\frac{\pal A}{\pal u^{(t)}(x)}
-\frac{\pal X(u(x); u_x(x),\dots)}{\pal u^{(t)}(x)} \pal_x^t A\right.\nn\\
&&\quad \quad \qquad\left. -\frac{\pal
X(u(y);u_y(y), \dots)}{\pal u^{(t)}(y)} \pal_y^t A\right).
\label{lie-der}
\end{eqnarray}
A bivector $\omega$ of the form (\ref{bi-vec}) defines a Poisson bracket
\begin{equation}
\{u(x),u(y)\}=
\sum_{k\ge 0} A_k(u(x);u_x(x),\dots) \delta^{(k)}(x-y)
\end{equation}
if and only if it satisfies the condition $[\omega, \omega]=0$.

Now let us denote by $\omega_1, \omega_2$ the bivectors corresponding to the two Poisson
brackets given in (\ref{z2}), and define the differentials $d_1, d_2:\ \Lambda_{loc}^k\to \Lambda_{loc}^{k+1}$ by
\begin{equation}
d_i \al=[\omega_i, \al],\ \al\in \Lambda_{loc}^k, \quad i=1,2;\, k\ge 0.
\end{equation}
The cohomologies of the complexes $(\Lambda_{loc}, d_i),\ i=1,2$ are called the Poisson cohomologies \cite{lichn}
of the hamiltonian structures $\omega_1, \omega_2$, and are denoted by $H^*(\omega_1), H^*(\omega_2)$.
Here $\Lambda_{loc}=\Lambda_{loc}^{0}\oplus\Lambda_{loc}^{1}\oplus\dots$.
The triviality of these Poisson cohomologies is proved in \cite{magri2,get} (see also \cite{DZ1}
for an alternative proof
of the triviality of the first two Poisson cohomologies).

\begin{prp}\label{mprp}
If the vector fields $\xi, \eta\in \Lambda^1_{loc}$ satisfy the
relation
\begin{equation}
d_1 \xi=d_2 \eta
\end{equation}
then there exist functionals $I, J$ of the form
\begin{equation}
I=\int f(u,u_x,\dots,u^{(m_1)}) dx,\quad J=\int g(u,u_x,\dots,u^{(m_2)}) dx
\end{equation}
such that
\begin{equation}
\xi=d_1 I-d_2 J.
\end{equation}
Here the densities are smooth functions of their arguments.
\end{prp}
Before we give the proof of the above proposition, let us first employ it to prove the Main Theorem.

\begin{prfM}
Due to the triviality of the Poisson cohomologies
$H^*(\omega_i),i=1,2$, the bihamiltonian
structure (\ref{z1}) can be transformed to the following form by a Miura transformation:
\begin{eqnarray}
\{u(x),u(y)\}_1&=&\vp(u(x))\,\delta'(x-y)+\frac12\p_x(\vp(u(x))\,\delta(x-y),\nn\\
\{u(x),u(y)\}_2&=&u(x)\,\vp(u(x))\,\delta'(x-y)+\frac12\p_x(u(x)\vp(u(x))\,\delta(x-y)
+\sum_{k \ge 1}\ve^k Q_k.\nn
\end{eqnarray}
Here $Q_k$ have the expression
$$
Q_k=\sum_{l=0}^{k+1}Q_{k,l}(u;u_x,\dots,u^{(k+1-l)})\,\delta^{(l)}(x-y).
$$
The compatibility of the above two Poisson brackets implies that
$$d_1 Q_1=0,\ d_2 Q_1=0$$
By using again the triviality of the Poisson cohomologies
$H^*(\omega_i),i=1,2$, we can find two local vector fields $\xi$ and $\eta$ such that
$$
Q_1=d_1 \xi,\ Q_1=d_2 \eta.
$$
It then follows from Proposition \ref{mprp} the existence of two local functionals $I, J$
satisfying $\xi=d_1I-d_2J$, this leads to the expression $Q_1=d_2d_1I$. So the term $\ve Q_1$ that appears in the
second Poisson bracket can be absorbed by the quasi-Miura transformation
$$
u\mapsto u-\ve\, d_1 I.
$$
After this quasi-Miura transformation, the first Poisson bracket is converted to the form
$$
\{u(x),u(y)\}_1=\vp(u(x))\,\delta'(x-y)+\frac12\p_x(\vp(u(x))\,\delta(x-y)+
\sum_{k \ge 2}\ve^k P_k.
$$
Here $P_k$ have similar expressions as that of $Q_k$. Note that the correction caused by the
above quasi-Miura transformation starts from the $\ve^2$ term.

By repeating the above procedure, we see that the deformation part
of (\ref{z1}) can be absorbed by a series of quasi-Miura transformations, so it is quasitrivial,
and the Main Theorem is proved.
\end{prfM}

Now let us proceed to prove the Proposition \ref{mprp}. To this end we first need to prove some lemmas.
We will use below the notation
\begin{equation}
u_k:=u^{(k)}=\p_x^k u(x),\quad k\ge 0.
\end{equation}
\begin{lem}\label{lem1}
If a local vector field $\xi$ has component of the form
\begin{equation}\label{form}
X=F(u,u_x,\dots,u_{[\frac{N}2]})\,u_N+Q(u,u_x,\dots,u_{N-1}),\
N \in \mathbb{N}
\end{equation}
then there exist two local functionals $I,J$
such that the component of the vector field $\xi-(d_1 I-d_2 J)$ depend at most on $u, u_x,\dots, u_{N-1}$.
\end{lem}

\begin{prf}
For any two local functionals $I, J$ of the form
$$
I=\int G(u,\dots,u_M) dx,\, J=\int H(u,\dots,u_M) dx
$$
we denote by $\tilde{\xi}$ the local vector field $d_1I-d_2J$ and by ${\tilde X}$ its component.
By a straightforward computation we have
\begin{eqnarray*}
(-1)^M\frac{\p \tilde{X}}{\p u_{2M+1}}&=&
\vp(u)\left(\frac{\p^2 G}{\p u_M \p u_M}-u\frac{\p^2 H}{\p u_M \p u_M}\right)\\
(-1)^M\frac{\p \tilde{X}}{\p u_{2M}}&=&
\left(M+\frac12\right)\vp(u)\,u_x\frac{\p^2 H}{\p u_M \p u_M}\\
&&+\left((M+1)\vp(u)\p_x+\frac12\vp'(u)u_x\right)\left(\frac{\p^2 G}{\p u_M \p
u_M} - u\frac{\p^2 H}{\p u_M \p u_M}\right)
\end{eqnarray*}
So when $N=2M+1$, the choice of $G, H$ that is given by
$$
G=\p^{-2}_{u_M}\left((-1)^M \frac{F(u,\dots,u_M)}{\vp(u)}\right),\ H=0
$$
yields the two local functionals $I, J$ that meet the requirement of the lemma. Similarly,
when $N=2M$, we can choose
$$
G=u H,\ H=\p^{-2}_{u_M}\left((-1)^M \frac{F(u,\dots,u_M)}{(M+1/2)\vp(u)u_x}\right).
$$
The lemma is proved.
\end{prf}

Let $\xi, \eta$ be two local vector fields with components $X, Y$ respectively of the form
\begin{equation}\label{xy-d}
X=X(u,\dots,u_N),\ Y=Y(u,\dots,u_N).
\end{equation}
Consider the bivector $d_1 \xi-d_2 \eta$, its component $Z$ can be written in the form
$$Z=\sum_{p \ge 0}Z_p \dl^{(p)}(x-y).$$
We denote by $Z_{p,s}$ the derivatives of $Z_p$ w.r.t $u_s$.
The main idea of our proof of Proposition \ref{mprp} is to reduce $X,Y$ to the form of
(\ref{form}) by using the vanishing of $Z_{p,s}$ when $d_1 \xi=d_2\eta $ holds true.

\begin{lem}\label{lem2}
If $Z_{0,2N+1}=0,\, Z_{0,2N}=0$, then the components of $\xi, \eta$ must take the following
form:
\begin{eqnarray}
&&X=(u G(u,\dots,u_{N-1})+F(u,\dots,u_{N-1}))u_N+Q(u,\dots,u_{N-1}),\nn\\
&&Y=G(u,\dots,u_{N-1})u_N+R(u,\dots,u_{N-1}).\label{xy-f}
\end{eqnarray}
\end{lem}
\begin{prf}
By a straightforward computation we have
$$
0=Z_{0,2N+1}=(-1)^{N+1}\vp(u)\left(\frac{\p^2 X}{\p u_N \p u_N}-u\frac{\p^2 Y}{\p u_N \p u_N}\right).
$$
So $X$ must take the form
$$
X=u Y+F(u,\dots,u_{N-1})u_N+{\tilde Q}(u,\dots,u_{N-1}).
$$
After substituting this expression
of $X$ into the expression of $Z$ we obtain
$$
0=Z_{0,2N}=(-1)^{N+1}\left(N+\frac12\right)\vp(u)u_x\frac{\p^2 Y}{\p u_N \p u_N}.
$$
It follows that $Y$ must have the expression
$$
Y=G(u,\dots,u_{N-1})u_N+R(u,\dots,u_{N-1}).
$$
The lemma is proved.
\end{prf}

\begin{lem}\label{lem3} Assume that the components $X, Y$ of the local vector fields $\xi,\eta$
have the form (\ref{xy-f}) and $N=2M+1$. Then
for any $m=1,2,\dots,M$ the following identity holds true:
\begin{equation}\label{ident}
\sum_{p=0}^m(-1)^{m-p}{N-p \choose m-p}Z_{p,2N+1-m-p}=\vp(u)\frac{\p F}{\p u_{N-m}}.
\end{equation}
\end{lem}
\begin{prf}
By our definition $Z$ is the component of the bivector $d_1\xi- d_2 \eta$, it has the explicit form
\begin{eqnarray*}
Z&=&\sum_{s \ge 0} (-1)^s \left(u\vp(u)\p_x^{s+1}(Y_s\dl)+
\frac12\p_x(u\vp(u))\p_x^{s}(Y_s\dl)\right)\\
&&-\sum_{s \ge 0} (-1)^s \left(\vp(u)\p_x^{s+1}(X_s\dl)+
\frac12\p_x(\vp(u))\p_x^{s}(X_s\dl)\right)+\cdots
\end{eqnarray*}
Here $\dl=\dl(x-y)$,\,$Y_s=\frac{\p Y}{\p u_s}$,\,$X_s=\frac{\p X}{\p u_s}$. There are some terms that are
omitted in the above expression, these terms do not affect the identity (\ref{ident}), so we can also omit them
in our calculations below.
Then $Z_p$ reads
\begin{eqnarray*}
Z_p&=&\sum_{s \ge 0} (-1)^s \left(u\vp(u){s+1\choose p}\p_x^{s+1-p}(Y_s)+
\frac12\p_x(u\vp(u)){s \choose p}\p_x^{s-p}(Y_s)\right)\\
&&-\sum_{s \ge 0} (-1)^s \left(\vp(u){s+1\choose p}\p_x^{s+1-p}(X_s)+
\frac12\p_x(\vp(u)){s\choose p}\p_x^{s-p}(X_s)\right).
\end{eqnarray*}
Denote by $l.h.s.$ the left hand side of the identity (\ref{ident}), then we obtain
\begin{eqnarray*}
l.h.s.&=&\sum_{p\ge0,s\ge0}(-1)^{m-p+s}{N-p \choose m-p}\\
&&\quad\left(u\vp(u){s+1\choose p}\sum_{t\ge0}{s+1-p\choose t}\p_x^t Y_{s,2N-m-s+t}\right.\\
&&\quad+\frac12\p_x(u\vp(u)){s\choose p}\sum_{t\ge0}{s-p\choose t}\p_x^t Y_{s,2N+1-m-s+t}\\
&&\quad-\vp(u){s+1\choose p}\sum_{t\ge0}{s+1-p\choose t}\p_x^t X_{s,2N-m-s+t}\\
&&\quad-\left.\frac12\p_x(\vp(u)){s\choose p}\sum_{t\ge0}{s-p\choose t}\p_x^t X_{s,2N+1-m-s+t}\right).
\end{eqnarray*}
Here we used the commutation relations
$$
\frac{\p}{\p u_l}\p_x^m=\sum_{t\ge0}{m\choose t}\frac{\p}{\p u_{l-m+t}}.
$$
By using the identity
$$
\sum_{p\ge0}(-1)^p{N-p \choose m-p}{s\choose p}{s-p\choose t}={s\choose t}{N-s+t\choose m}
$$
and by changing the order of summation, we can rewrite $l.h.s.$ as follows
\begin{eqnarray*}
&&(-1)^m l.h.s=\sum_{t\ge0}
\left(u\vp(u)\p_x^t\left(\sum_{s\ge0}(-1)^s{s+1\choose t}{N-s+t-1\choose m}Y_{s,2N-m-s+t}\right)\right.\\
&&\qquad+\frac12\p_x(u\vp(u))\p_x^t\left(\sum_{s\ge0}(-1)^s{s\choose t}{N-s+t\choose m}Y_{s,2N+1-m-s+t}\right)\\
&&\qquad-\vp(u)\p_x^t\left(\sum_{s\ge0}(-1)^s{s+1\choose t}{N-s+t-1\choose m}X_{s,2N-m-s+t}\right)\\
&&\qquad\left.
-\frac12\p_x(\vp(u))\p_x^t\left(\sum_{s\ge0}(-1)^s{s\choose t}{N-s+t\choose m}X_{s,2N+1-m-s+t}\right)\right).
\end{eqnarray*}
Now we substitute the expression (\ref{xy-f}) of $X,Y$ into the right hand side of the above formula.
It's easy to see that all terms in the above summation vanish
except the terms with $s=N,t=0$, so the above formula can be simplified to
\begin{eqnarray*}
l.h.s&=&(-1)^m\left((-1)^{N+m}u\vp(u)Y_{N,N-m}-(-1)^{N+m}\vp(u)X_{N,N-m}\right)\\
&=&\vp(u)\frac{\p F}{\p u_{N-m}}.
\end{eqnarray*}
The lemma is proved.
\end{prf}

For two local vector fields $\xi, \eta$ whose components have the form (\ref{xy-d}) with $N=2M+1$, the above lemmas show
that $X, Y$ must have the form (\ref{xy-f}) and $F$ is independent of $u_{M+1},\dots,u_{N}$.
From the proof of Lemma \ref{lem1}, we know that we can find a local functional $I$ such that the component
of $\xi-d_1I$ has the expression of (\ref{xy-f}) with vanishing $F$.
On the other hand, when $N$ is even, because
$Z_{N+1}=-2\vp(u)F$, we immediately have $F=0$. So for both cases with even $N$ and odd $N$, we can always
modify, if necessary, the vector filed $\xi$ by subtracting a vector field of the form $d_1 I$ such that
the components
$X,Y$ of the local vector fields $\xi, \eta$ can be expressed as
\begin{eqnarray}
&&X=u G(u,\dots,u_{N-1})u_N+Q(u,\dots,u_{N-1}),\nn \\
&&Y=G(u,\dots,u_{N-1})u_N+R(u,\dots,u_{N-1}).\label{xy-g}
\end{eqnarray}

\begin{lem}\label{lem4}
Suppose the function $G$ that appears in (\ref{xy-g}) depends only on $u, u_x,\dots,$\newline $ u_{N-m}$ with
$1\le m\le [\frac{N-1}2]$, then we have the following identity:
\begin{equation}
\sum_{p=0}^m(-1)^{m-p}{N-p \choose m-p}Z_{p,2N-m-p}=(-1)^{N+1}\left(N-m+\frac12\right)\vp(u)u_x\frac{\p G}{\p u_{N-m}}.
\end{equation}
\end{lem}

\begin{prf}
The proof of the above identity is similar to that of (\ref{ident}). We only need to note that
in the last step of the proof the terms with $(t,s)=(1,N-m),(0,N),(0,N-m-1),(0,N-m)$ do not
vanish. We omit the details of the calculations here. The lemma is proved.
\end{prf}

From the last lemma we see that the function $G$ depends at most on
$u, u_x,\dots,$\newline $ u_{[\frac{N}2]}$. So it follows from Lemma \ref{lem1} that when
$N=2 M+1$ the components of the local vector fields ${\hat \xi}=\xi+d_2 I_1, {\hat \eta}=\eta-d_1 I_1$ with
$$
I_1=(-1)^M \int \p^{-2}_{u_M}\frac{G(u,u_x,\dots, u_M)}{\vp(u)} dx
$$
depend at most on $u, u_x,\dots, u_{N-1}$, while the relation $d_1{\hat \xi}=d_2 {\hat \eta}$ still holds true.
When $N=2M$, we can consider the following modified local vector fields
$$
{\hat \xi}=\xi-(d_1 I_1-d_2 J_1),\quad {\hat \eta}=\eta-(d_1 J_1-d_2 J_2).
$$
Here
$$
I_1=\int u(x)^2 P dx,\quad J_1=\int u(x) P dx,\quad J_2=\int P dx
$$
with
$$
P=(-1)^M \p^{-2}_{u_M}\frac{G(u,u_x,\dots, u_M)}{(M+1/2) \vp(u)\,u_x}.
$$
The components of these modified local vector fields depend also at most on
$u, u_x,\dots, u_{N-1}$, and the relation $d_1{\hat \xi}=d_2 {\hat \eta}$ is unchanged.
Now the proof of Proposition \ref{mprp} is easily obtained by induction on the integer
$N$ that appears in the expression (\ref{xy-d}) of the components of the vector fileds $\xi, \eta$.

\section{Conclusion}
In the proof of the quasitriviality of the deformed bihamiltonian structure (\ref{z1}) we
do not use the fact that $A_{m,l;a}$ are differential polynomials, so the above result of quasitriviality
can be extended to the class of deformed bihamiltonian structures of the form (\ref{z1}) without the assumption of
polynomiality.

The above method can be employed to study the quasitriviality of
deformations of the general multi-component bihamiltonian structures of the form
\begin{eqnarray}\label{z3}
&&\{u^i(x),u^j(y)\}_a=g^{ij}_a(u(x))\delta'(x-y)+\Gamma^{ij}_{k;a}(u(x)) u^k_x \delta(x-y),\nn\\
&&\quad
i,j=1,\dots, n,\ a=1,2.\nn
\end{eqnarray}
Here $(g^{ij}_1), (g^{ij}_2)$ form a flat pencil of metrics\cite{Du1} on a manifold $M$, and $\Gamma^{ij}_{k;1},
\Gamma^{ij}_{k;2}$ are the corresponding contravariant
coefficients of the Levi-Civita connections of these flat metrics. It is conjectured in \cite{LZ} that
for such a bihamiltonian structure all its deformations of the form that is similar to (\ref{z1})
are quasitrivial. Under the assumption of validity of this conjecture, it is proved in \cite{LZ}
that any equivalence class of such deformations is uniquely determined by $n$ functions
of one variable. We will return to the proof of this conjecture in a subsequent publication.

\vskip 0.2truecm \noindent{\bf Acknowledgments.} The authors are
grateful to Boris Dubrovin for helpful suggestions and comments.
The researches of Y.Z. were partially
supported by the Chinese National Science Fund for Distinguished
Young Scholars grant No.10025101 and the Special Funds of Chinese
Major Basic Research Project ``Nonlinear Sciences''.

\end{document}